\documentclass[12pt]{article}
\usepackage[utf8]{inputenc}
\usepackage{parskip}
\usepackage{amsmath,amsfonts,amssymb,amsthm,mathtools}
\usepackage{commath}
\usepackage{fullpage}
\usepackage{hyperref}
\usepackage[numeric]{amsrefs}
\newcommand{\Z}{\mathbb Z}
\newcommand{\R}{\mathbb R}
\newcommand{\N}{\mathbb N}
\newcommand{\Q}{\mathbb Q}

\newcommand{\inR}{\int_{-\infty}^{\infty}}
\newcommand{\inP}{\int_{0}^{\infty}}

\newcommand{\intres}{\int_{e^{-\frac{2}{T}}}^{e^{\frac{2}{T}}}}
\newcommand{\inIm}{\int_{-i\infty}^{i\infty}}
\newcommand{\sumn}{\sum_{n\geq 1}}

\newcommand{\dx}{\,dx}
\newcommand{\dt}{\,dt}
\newcommand{\ds}{\,ds}

\newcommand{\xdx}{\,\frac{dx}{x}}

\newcommand{\SL}{\text{SL}}
\newcommand{\GL}{\text{GL}}
\newcommand{\real}{\text{Re}\,}

\newcommand{\sgn}{\text{sgn}\,}

\newtheoremstyle{theoremstyle}{}{}{\itshape}{}{\scshape}{.}{ }{\textbf{#1\ #2}}  
\theoremstyle{plain}

\numberwithin{equation}{section}

\newtheorem*{theorem*}{Theorem}
\newtheorem{theorem}[equation]{Theorem}
\newtheorem{lem}[equation]{Lemma}
\newtheorem*{lem*}{Lemma}

\newtheorem*{prop*}{Proposition}
\newtheorem{cor}[equation]{Corollary}

\theoremstyle{definition}
\newtheorem{defin}[equation]{Definition}

\numberwithin{equation}{section}

\begin{document}
\title{Infinitely many zeros of additively twisted $L$-functions on the critical line}
\author{Doyon Kim}
\maketitle
\begin{abstract}
For $f$ a cuspidal modular form for the group $\Gamma_0(N)$ of integral or half-integral weight, $N$ a multiple of $4$ in case the weight is half-integral, we study the zeros of the $L$-function attached to $f$ twisted by an additive character $e^{2\pi i n \frac{p}{q}}$ with $\frac{p}{q}\in \Q$. We prove that for certain $f$ and $\frac{p}{q}\in \Q$, the additively twisted $L$-function has infinitely many zeros on the critical line. We develop a variant of the Hardy-Littlewood method which uses distributions to prove the result.
\end{abstract}
\section{Introduction}
\subsection{The Hardy-Littlewood method}
Hardy and Littlewood \cite{hardyzeta} proved that the Riemann $\zeta$-function vanishes infinitely many times on the critical line by constructing a continuous real-valued function $Z(t)$ defined on the real line whose zeros correspond to the critical zeros of the $\zeta$-function, and showing that $Z(t)$ changes sign infinitely often. They used
\begin{equation}\label{Ztzeta}
    Z(t)=\frac{H(1/2+it)}{\abs{H(1/2+it)}}\zeta(1/2+it),
\end{equation}
where
\begin{equation}\label{Htzeta}
    H(s)=\frac{1}{2}s(1-s)\pi^{-s/2}\Gamma(s/2).
\end{equation}
The real-valuedness of $Z(t)$ follows from the functional equation 
\begin{equation}\label{fneqzeta}
H(s)\zeta(s)=H(1-s)\zeta(1-s).\end{equation} A sign change of $Z(t)$ in an interval $(a,b)$ can be detected by comparing the size of the two integrals
\[
    \int_{a}^{b} Z(t)\dt \quad \text{and} \quad  \int_{a}^{b} \abs{Z(t)}\dt.
\]
The two integrals have different magnitudes if and only if $Z(t)$ changes sign somewhere in the interval. If $Z(t)$ has only finitely many zeros on the critical line, then $Z(t)$ does not change sign in the interval $[T,2T]$ for sufficiently large $T$, hence
\begin{equation}\label{intT2T}
     \abs{\int_{T}^{2T} Z(t)\dt} =\int_{T}^{2T} \abs{Z(t)}\dt.
\end{equation}
Hardy and Littlewood showed that the integral on the right hand side grows at least as fast as a constant multiple of $T$ as $T\to\infty$ and the integral on the left hand side is bounded by a fractional power of $T$, therefore $Z(t)$ has infinitely many zeros on the real line. \par
This method is called the Hardy-Littlewood method, and this method has been used for different classes of $L$-functions as well. Epstein, Hafner and Sarnak \cite{epsteinmaass}, \cite{hafnermaass} proved that the $L$-functions of Maass cusp forms have infinitely many zeros on the critical line, and Chandrasekharan and Narasimhan \cite{chanquadratic} showed that the zeta functions of ideal classes in quadratic fields by have infinitely many zeros on the critical line. In both cases the real-valued functions $Z(t)$ of are of the form
\begin{equation}
    Z(t)=\chi_L(1/2+it)^{-1/2}L(1/2+it),
\end{equation}
where $\chi_L(s)$ is the holomorphic function satisfying the functional equation
\begin{equation}
    L(s)=\chi_L(s)L(1-s).
\end{equation}
 The most technical part of the Hardy-Littlewood method is obtaining the upper bound of the integral
\[
    \int_{T}^{2T} Z(t)\dt,
\]
which requires the cancellation of exponential sums
\begin{equation}\label{expsum}
    \sum_{m\leq T} a_m e^{2\pi m x}=o(T), \quad \text{as}\quad T\to\infty,
\end{equation}
for any value of $x$. \par 
Wilton \cite{wiltonadditive} used a variant of the Hardy-Littlewood method to prove that if $p/q\in\Q$ is a rational number such that $p^2\equiv 1\pmod q$, then the additively twisted $L$-function of modular discriminant $\Delta(z)$
\begin{equation}
    \rho(s;p/q)=\sumn \frac{\tau(n)}{n^s}e^{2\pi i n p/q},
\end{equation}
where $\tau(n)$ denotes the $n$-th Fourier coefficient of $\Delta(z)$, has infinitely many zeros on the line $\real s=6$. In his proof, Wilton used the function
\begin{equation}
    Z_{\alpha,n}(t)=t^n R(t) e^{\alpha t},
\end{equation}
where
\begin{equation}
R(t)=(2\pi/q)^{-6-it}\Gamma(6+it)L_{p/q}(6+it).    
\end{equation}
It follows from a functional equation of $L_{p/q}(s,\Delta)$ that $Z_{\alpha,n}(t)$ is real-valued for $t\in\R$. Wilton showed that
\begin{equation}
    \lim_{\alpha\to\pi/2^-}\inR t^n R(t) e^{\alpha t}\dt=0,
\end{equation}
for any positive integer $n$, using the limiting behavior of higher order derivatives of the modular discriminant near the cusp $p/q$. From this he deduced the fact that $L_{p/q}(s,\Delta)$ has infinitely many zeros on the critical line. \par
Recently, Meher, Pujahari and Kotyada \cite{meherhalfcusp1} proved that $L$-functions of half-integral weight cusp forms of level $4$ have infinitely many zeros on the critical line, and Meher, Pujahari and Shankhadhar \cite{meherhalfcusp2} generalized the result to cusp forms of levels $4N^2$. In \cite{meherhalfcusp1}, the proof involves establishing a lower bound on an exponential sum of the form (\ref{expsum}) and in \cite{meherhalfcusp2} the proof uses Wilton's variant. \par
\subsection{Main results}
In this paper, we use a distributional variant of the Hardy-Littlewood method to generalize the result of Wilton \cite{wiltonadditive} to cusp forms of integral or half-integral weight and higher level. The main advantage of using distributions is that it does not require any bound of the exponential sums, making it easily applicable to different classes of cusp forms, including those of half-integral weights. Instead of the regularity of exponential sums, we use the notion of a distribution vanishing to infinite order. \par 
For a complex number $z$, we put $e(z)=e^{2\pi i z}$. For $\nu\in \frac{1}{2}\Z$ and $N$ a positive integer with a condition that $4|N$ if $\nu\notin \Z$, we denote $S_\nu\left(\Gamma_0(N)\right)$ the space of weight $\nu$ cusp forms with respect to the congruence subgroup $\Gamma_0(N)$. Since the matrix $\begin{pmatrix} 1 & 1 \\ 0 & 1 \end{pmatrix}$ is an element of $\Gamma_0(N)$, a cusp form $f(z)$ in $\Gamma_0(N)$ has a Fourier expansion of the form
\begin{equation}f(z)=\sumn a_n n^{\frac{\nu-1}{2}}e(nz).\end{equation}
Let $L(s)=L(s,f)$ denote the $L$-function of $f(z)$,
\begin{equation}
    L(s)=L(s,f)=\sumn a_n n^{-s}.
\end{equation}
Here, we normalize the Fourier coefficients so the critical line is the vertical line $\real s=\frac{1}{2}$ for any weight $\nu$. For a rational number $\frac{p}{q}\in\Q$, we let $L_{p/q}(s)=L_{p/q}(s,f)$ denote the additively twisted $L$-function
\begin{equation}
    L_{p/q}(s)=\sumn a_n e(np/q) n^{-s}.
\end{equation}
\begin{theorem}\label{integralthm}
Let $k$ be an even positive integer, $N$ be a positive integer, and $\frac{p}{q}$ a rational number which is $\Gamma_0(N)$-equivalent to $i\infty$ such that $p^2\equiv 1 \pmod q$. If $f(z)\in S_k\left(\Gamma_0(N)\right)$ has a Fourier expansion
    \begin{equation}f(z)=\sumn a_n n^{\frac{k-1}{2}}e(nz),\end{equation}
with $a_1\neq 0$ and $a_j\in\R$ for all $j$, then $L_{p/q}(s,f)$ has infinitely many zeros on the critical line.
\end{theorem}
\begin{theorem}\label{halfintegralthm} Let $k$ be a positive integer, $N$ a positive integer divisible by $4$, and $\frac{p}{q}$ a rational number which is $\Gamma_0(N)$-equivalent to $i\infty$ such that $p^2\equiv 1 \pmod q$.
If $g(z)\in S_{k+\frac{1}{2}}\left(\Gamma_0(N)\right)$ has a Fourier expansion
    \begin{equation}g(z)=\sumn b_n n^{\frac{k}{2}-\frac{1}{4}}e(nz),\end{equation}
with $b_1\neq 0$ and $b_j\in\R$ for all $j$, then $L_{p/q}(s,g)$ has infinitely many zeros on the critical line.
\end{theorem}
\section{Automorphic distributions attached to cusp forms}
\subsection{Automorphic distributions on \texorpdfstring{$SL(2,\Z)$}{}}
Let $f(z)$ be a cusp form of weight $k$ for $SL(2,\Z)$ on the upper half plane $\mathbb{H}$ with Fourier expansion 
\begin{equation}\label{ffourier}
f(z)=\sum_{n\geq 1}a_n n^{\frac{k-1}{2}}e(nz).\end{equation}
Then $f(z)$ satisfies the automorphy condition 
\begin{equation}\label{fautomorphy}
f\left(\frac{az+b}{cz+d}\right)=(cz+d)^{k}f(z), \quad \text{for all}\quad \begin{pmatrix} a & b \\ c& d \end{pmatrix}\in \SL(2,\Z). \end{equation}
The cusp form $f(z)$ has distribution boundary values. That is, the limit
\begin{equation}\label{taulimdef}
\tau(x)=\tau_f(x)=\lim_{y\to 0^{+}} f(x+iy).\end{equation}
converges in the strong distribution topology.
As a consequence of the limit formula, $\tau$ inherits the automorphicity from $f(z)$, so
\begin{equation}\label{tauautomorphy}
\tau\left(\frac{ax+b}{cx+d}\right)=(cx+d)^{k}\tau(x),\quad \text{for all}\quad \begin{pmatrix} a & b \\ c& d \end{pmatrix}\in \SL(2,\Z).\end{equation}
We call such distribution $\tau$ the automorphic distribution attached to $f(z)$.
The automorphy condition (\ref{tauautomorphy}) with $a=b=d=1$ and $c=0$ implies \[\tau(x+1)=\tau(x),\] so $\tau$ can be written as a Fourier series. The limit (\ref{taulimdef}) can be taken from the Fourier series expansion (\ref{ffourier})  term by term, and we have
\begin{equation}\label{taufourier}
\tau(x)=\sum_{n\geq 1}a_n n^{\frac{k-1}{2}}e(nx).\end{equation}
Since $\tau$ has no constant term, $\tau$ can be written as the $j$-th derivative of a continuous, periodic function for every sufficiently large integer $j$, 
\begin{equation}\label{taujth}
\tau^{(-j)}(x)=\sumn (2\pi in)^{-j}a_n n^{\frac{k-1}{2}}e(nx).\end{equation}
The distribution $\tau$ acts on any Schwartz function $\psi(x)\in\mathcal{S}(\R)$ by integration by parts,
\begin{equation}\label{ibp}
    \inR \tau(x)\psi(x)\dx=(-1)^j\inR \tau^{(-j)}(x)\frac{d^j}{dx^j}\psi(x)\dx.
\end{equation}
Furthermore, $\tau$ can be integrated against the Mellin kernel $\abs{x}^{s-1}$, even though it is not a Schwartz function. This enables us to define Mellin transform of $\tau$. The Mellin transform of periodic distributions will be discussed in Section 2.3.

\subsection{Order of vanishing of distributions}

In this section, we define the notion of vanishing of a distribution $\sigma(x)$ to order $k\geq 0$ and present how the notion can be used to derive bounds of certain integrals. The following definition is from Definition 2.4 and Lemma 3.1 of \cite{inforder}
\begin{defin}\label{vanishingorder}
A distribution $\sigma(x)\in C^{-\infty}(\R)$ vanishes to order $m\geq 0$ at a point $p\in\R$ if there exists an open interval $I_p$ containing $p$ and locally bounded functions $h_j\in L^{\infty}_{\text{loc}}(I_p)$ indexed by $0\leq j\leq N$ such that 
\begin{equation}\label{distsum}
    \sigma(x)=\sum_{j=0}^N \frac{d^j}{dx^j}\left((x-p)^{m+j}h_j(x)\right)
\end{equation}
as an identity between distributions on $I_p$. The distribution $\sigma$ vanishes to infinite order at $p$ if it vanishes to order $m$ for every $m\geq 0$.
\end{defin}
In case $\sigma(x)$ vanishes to order $m+1$ at $p$, we have a stronger result.
\begin{lem}[Lemma 3.1 of \cite{inforder}]
If $\sigma(x)\in C^{-\infty}(\R)$ vanishes to order $m+1$ at a point $p\in\R$, then we can write 
\begin{equation}
    \sigma(x)=\sum_{j=0}^N \frac{d^j}{dx^j}\left((x-p)^{m+j}h_j\right)
\end{equation}
with $h_j\in C(I_p)$.
\end{lem}
We define the order of vanishing of $\sigma(x)$ at $\infty$ as the order of vanishing of $\sigma(1/x)$ at $0$. Using the definition above, we prove the following lemma.
\begin{lem}\label{tauinforderlem} The automorphic distribution $\tau$ attached to a cusp form  $f(z)$ for $\SL(2,\Z)$ vanishes to infinite order at every cusp.
\end{lem}
\begin{proof}
The automorphicity of $\tau(x)$
combined with the chain rule for the change of variables $x$ to $(ax+b)/(cx+d)$ implies
\begin{equation}\label{tauchainrule}
\tau(x)=(cx+d)^{-k}\left((cx+d)^2\frac{d}{dx}\right)^l\tau^{(-l)}\left(\frac{ax+b}{cx+d}\right)\end{equation}
for every sufficiently large $l\in\N$. Expanding the differential operator, we obtain the expression
\begin{equation}\label{tauchain2}
\tau(x)=\sum_{j=1}^{l} q_j (cx+d)^{l-k+j}\frac{d^j}{dx^j}\left(\tau^{(-l)}\left(\frac{ax+b}{cx+d}\right)\right).\end{equation}
Since $\tau^{(-l)}(x)$ is a continuous periodic function, 
$\tau^{(-l)}\left((ax+b)/(cx+d)\right)$ is bounded near $-d/c$. Therefore $\tau^{(-l)}\left((ax+b)/(cx+d)\right)$ has a locally continuous antiderivative. This gives the expression
\begin{equation}\label{tauchain3}
\tau(x)=\sum_{j=1}^{l} q_j (cx+d)^{l-k+j}\frac{d^{j+1}}{dx^{j+1}}H(x)\end{equation}
where $H(x)$ is continuous near $-d/c$. Moving the factor $(cx+d)^{l-k+j}$ across the differential operator, we can write $\tau$ as
\begin{equation}\label{tausum2}
\tau(x)=\sum_{j=0}^{l+1} \frac{d^j}{dx^j}\left(x+\frac{d}{c}\right)^{l-k+j-1} h_j(x),\end{equation}
where $h_j(x)$ are continuous functions on the neighborhood of $-d/c$. \par
Let $m$ be a positive integer. The presentation (\ref{tausum2}) with $l=m+k+1$ gives
\begin{equation}\label{tausum}
\tau(x)=\sum_{j=0}^{m+k+2} \frac{d^j}{dx^j}\left(x+\frac{d}{c}\right)^{m+j} h_j(x),\end{equation}
where $h_j(x)$ are continuous on the neighborhood of $-d/c$. Therefore, the distribution $\tau(x)$ vanishes to order $m$ at $-d/c$. Since the choice of $m$ and the coprime pair $c$ and $d$ are arbitrary, we deduce that $\tau(x)$ vanishes to infinite order at every rational. Similarly, it follows from the equation
\begin{equation}\label{tauatinfty}
\tau\left(-\frac{1}{x}\right)=\left(x^2\frac{d}{dx}\right)^l\tau^{(-l)}\left(-\frac{1}{x}\right)\end{equation} 
that $\tau(x)$ vanishes to infinite order at $\infty$.
\end{proof}
The equation (\ref{tauatinfty}) holds for any periodic distribution without the constant term. Hence any periodic distribution without the constant term vanishes to infinite order at $\infty$.
\begin{lem}\label{vanishingboundlem}
Suppose that $\sigma$ is a distribution which vanishes to order $m\geq 0$ at $p$. Let $\left\{\varphi_T\right\}_{T>0}$ be a family of smooth functions such that $\varphi_T(x)$ is supported on $(p-c_1/T,p+c_1/T)$ for some $c_1>0$ and $\abs{\frac{d^j}{dx^j}\varphi_T(x)}\leq c_2 T^{j+1}$ for some $c_2>0$. Then
\begin{equation}\label{vanishingbound}
    \inR \sigma(x) \varphi_T(x)\dx\ll_m T^{-m}.
\end{equation}
as $T\to \infty$. In particular, if the distribution $\sigma$ vanishes to infinite order at $p$, then the integral above decays rapidly as $T\to\infty$, that is, the integral is $O(T^{-m})$ for any $m>0$.
\end{lem}
\begin{proof}
By Definition~\ref{vanishingorder}, there exists an $\epsilon>0$ such that we can write $\sigma$ as
\begin{equation}
    \sigma(x)=\sum_{j=0}^{N} \frac{d^j}{dx^j}(x-p)^{m+j}h_j(x) 
\end{equation}
where $h_j\in L_{\text{loc}}^{\infty} \left((p-\epsilon,p+\epsilon)\right)$. Take large $T>0$ so that $c_1/T<\epsilon$. We have
\begin{equation}
    \begin{aligned}
\inR \sigma(x)\varphi_T(x)\dx&=\int_{p-c_1/T}^{p+c_1/T} \left(\sum_{j=0}^{N} \frac{d^j}{dx^j}(x-p)^{m+j}h_j(x)\right)\varphi_T(x)\dx \\
&=\sum_{j=0}^{N}\int_{p-c_1/T}^{p+c_1/T} \left(\frac{d^j}{dx^j}(x-p)^{m+j}h_j(x)\right)\varphi_T(x)\dx \\
&=\sum_{j=0}^{N} (-1)^j \int_{p-c_1/T}^{p+c_1/T} (x-p)^{m+j}h_j(x) \frac{d^j}{dx^j}\varphi_T(x)\dx.
    \end{aligned}
\end{equation}
Observe that each integral in the summand is $O(T^{-1}\cdot T^{-m-j} \cdot T^{j+1})=O(T^{-m})$. The Lemma follows.
\end{proof}

\subsection{Mellin transform of periodic distributions}
For $\delta\in \Z/2\Z$, we define the signed Mellin transform of a distribution $\sigma$ as
\begin{equation}\label{Mellindef}
M_\delta\sigma(s)=\int_{-\infty}^{\infty} \sgn (x)^\delta\sigma(x)\abs{x}^{s-1}\dx.\end{equation}
Computing formally, without worrying about the convergence, the signed Mellin transform of a periodic distribution $\tau(x)=\sumn a_n e(nx)$ gives
\begin{equation}\label{Mellinformal}
    \begin{aligned}
 M_\delta\tau(s)&=\inR \sgn(x)^\delta\biggl(\sumn a_n e(nx)\biggr)\abs{x}^{s-1}\dx =\sumn a_n \inR\sgn(x)^\delta e(nx)\abs{x}^{s-1}\dx \\
 &=\sumn a_n n^{-s}\inR \sgn(x)^\delta e(x)\abs{x}^{s-1}\dx=G_{\delta}(s)\sumn a_n n^{-s},
    \end{aligned}
\end{equation}
where
\begin{equation}\label{distgamma}
\begin{aligned}
G_\delta(s)=2i^\delta (2\pi)^{-s}\Gamma(s)\cos\left(\frac{\pi(s-\delta)}{2}\right).
\end{aligned}
\end{equation}
The following Lemma, which is from \cite{inforder}, legitimizes the formal computation on a right half plane.
\begin{lem}[Lemma 3.38 of \cite{inforder}]\label{Mellinrightlem}
Let $\tau$ be a periodic distribution without the constant term written as a Fourier series $\tau(x)=\sumn a_n e(nx)$, such that its Fourier coefficients satisfy the bound $a_n=O(n^\alpha)$ for some $\alpha>0$. Then the signed Mellin transform of $\tau(x)$ is given by the formula 
\begin{equation}
    M_\delta\tau(s)=G_\delta(s)\sumn a_n n^{-s} \quad \text{for} \quad \real s\gg 1.
\end{equation}
\end{lem}
The following lemma extends the Mellin transform to the entire complex plane. The Lemma is from Section 3 of \cite{inforder}.
\begin{lem}\label{Mellinentirelem}
Let $\tau$ be as in Lemma~\ref{Mellinrightlem}, and further assume that $\tau$ vanishes to infinite order at $0$. Then the signed Mellin transform extends to an entire function.
\end{lem}
\begin{proof}
It is shown in Section 3 of \cite{inforder} that if a distribution $\sigma$ vanishes to order $k_0\geq 0$ and has a canonical extension to $\infty$ which vanishes there to order $k_{\infty}\geq 0$, then the signed Mellin transform given by (\ref{Mellindef}) is holomorphic on $-k_0<\real s<k_{\infty}$. Since a periodic distribution without constant term vanishes to infinite order at $\infty$, we have $k_0=k_{\infty}=\infty$.
\end{proof}

We state Parseval's identity for Mellin transforms of distributions. For a Schwartz function $\psi(x)$, we let $M\psi(s)$ without subscript $\delta\in \Z/2\Z$ denote the usual Mellin transform
\begin{equation}\label{Mellinschwartz}
M\psi(s)=\inP \psi(x)x^{s-1}\dx.\end{equation}
\begin{lem}[Lemma 4.25 of \cite{inforder}]\label{Parsevallem}
Let $k\geq 0$ and let $\sigma$ be a distribution that vanishes to order $k$ at $0$ and has an extension across $\infty$ vanishing to order $k_{\infty}\geq 1$ at $\infty$. Then the identity
\[\inR \psi(x)\sigma(x)\dx=\frac{1}{4\pi i}\int_{\real s=s_0}M_{\delta}\psi(s)M_{\delta}\sigma(1-s)\ds,\]
holds for all Schwartz functions $\psi$ satisfying $f(-x)=(-1)^\delta f(x)$ and $0<s_0<k+1$.
\end{lem}
\begin{cor}\label{Parsevalcor}
    Let $\sigma$ be a distribution that vanishes to infinite order at $0$ and has an extension across $\infty$ vanishing to infinite order at $\infty$. Let $\psi$ be a smooth function of compact support which vanishes near $0$. The identity
    \begin{equation}\label{Parsevalid}
    \inR \sigma(x) \psi(\abs{x}) \abs{x}^{k-1}\dx=\frac{1}{2\pi i}\inIm M_0\sigma(k-s)M\psi(s)\ds
\end{equation}
holds for all $k>0$.
\end{cor}
\begin{proof}
Let $\widetilde{\psi}(x)=\psi(\abs{x})$. Then $\widetilde{\psi}$ is an even Schwartz function and satisfies
\[M_0\widetilde{\psi}(s)=\inR \psi(\abs{x})\abs{x}^{s-1}\dx=2\inP \psi(x) x^{s-1}\dx=2M\psi(s).\]
It follows from Definition~\ref{vanishingorder} that as a distribution, $\sigma(x)\abs{x}^{k-1}$ vanishes to infinite order at $0$ and has an extension across $\infty$ vanishing to infinite order at $\infty$. Applying Lemma~\ref{Parsevallem} with $\widetilde{\psi}$ and $\sigma(x)\abs{x}^{k-1}$ in the place of $\psi$ and $\sigma$ yields (\ref{Parsevalid}).
\end{proof}

\subsection{Outline of the proof}\label{outline}
We end Section 2 with an outline of our proof of the main theorems. Let $\left\{u_T\right\}_{T>0}$ be a family of real-valued, positive and uniformly bounded functions which satisfies the bound
\begin{equation}\label{testfunction}
    u_T(x)\ll_N T^{N}\abs{x-2T^{3/2}}^{-N}
\end{equation}
for any $N>0$. In particular, $u_T(x)\ll_N T^{-N}$ outside of the interval $(T^{3/2},3T^{3/2})$. If $f(x)$ is a real-valued function which grows at most polynomially as $\abs{x}\to\infty$, then the following bound holds:
\begin{equation} \label{testbound1}
    \int_{\abs{x-2T^{3/2}}>T^{3/2}} \abs{f(x)}u_T(x)\dx\ll_N T^{-N}.
\end{equation}
Furthermore, if $f(x)$ does not vanish in the interval $(T^{3/2},3T^{3/2})$, then we have
\begin{equation}\label{testbound2}
    \inR \abs{f(x)} u_T(x)\dx-\inR f(x) u_T(x)\dx\ll_N T^{-N}.
\end{equation}
In Section 3 and 4, we construct a real-valued function $Z(x)$ whose zeros correspond to the critical zeros of an additively twisted $L$-function, using the Mellin transform of the associated automorphic distribution. In Section 5, we construct a family of test functions $u_T(x)$ parameterized by $T$ using bump functions, and establish the bounds (\ref{testbound1}) and (\ref{testbound2}). Theorem~\ref{integralthm} and Theorem~\ref{halfintegralthm} follow once we establish the bounds
\begin{equation}\label{bound1goal}
\inR Z(x) u_T(x)\dx= o(T)
\end{equation}
and
\begin{equation}\label{bound2goal}
\inR \abs{Z(x)} u_T(x)\dx\sim T,
\end{equation}
which together imply that $Z(x)$ has a zero in the interval $(T^{3/2},3T^{3/2})$. Since this is true for any sufficiently large $T$, it follows that $Z(x)$ has infinitely many zeros on the real line. We obtain the bound (\ref{bound1goal}) using Parseval's identity (\ref{Parsevalid}) and the argument presented in Lemma~\ref{vanishingboundlem}. The bound (\ref{bound2goal}) is established using the Mellin inversion formula.

\section{Modular forms of integral weight}
In this section, we consider automorphic distributions attached to cusp forms of weight $k\in\Z$ for the modular subgroup $\Gamma_0(N)$. Since $-I\in\Gamma_0(N)$,  $S_k\left(\Gamma_0(N)\right)=\{0\}$ if $k$ is odd. So throughout the section we assume that $k$ is even. Expand $f(z)$ in a Fourier series at the cusp $i\infty$ as
\begin{equation}\label{ffourierlevelN}
f(z)=\sum_{n\geq 1}a_n n^{\frac{k-1}{2}}e(nz).\end{equation}
As before, we define the automorphic distribution $\tau=\tau_f$ attached to $f$ to be \begin{equation}\label{taufourierlevelN}
\tau(x)=\lim_{y\to 0^+}f(x+iy)=\sum_{n\geq 1}a_n n^{\frac{k-1}{2}}e(nx).\end{equation}
The distribution $\tau$ satisfies the automorphy relation
\begin{equation}\label{tauautomorphylevelN}
\tau\left(\frac{ax+b}{cx+d}\right)=(cx+d)^{k}\tau(x),\quad \begin{pmatrix} a & b \\ c& d \end{pmatrix}\in \Gamma_0(N).\end{equation}

For $\alpha=\begin{pmatrix} a & b \\ c & d \end{pmatrix}\in\SL(2,\Z)$, let $j(\alpha,z)=cz+d$. Define slash operator $f|\alpha$ by
\begin{equation}\label{intslash}
(f|\alpha) (z)=j(\alpha,z)^{-k}f(\alpha z).\end{equation}
Then $f|\alpha$ is modular for the group $\alpha^{-1}\Gamma_0(N)\alpha$. Since $f(z)$ is a cusp form, $f|\alpha$ has a Fourier expansion 
\begin{equation}
(f|\alpha)(z)=\sumn b_n e(nz/h)\end{equation} for some $h\in\Z$ and has a corresponding automorphic distribution 
\begin{equation}(\tau|\alpha)(x)=\lim_{y\to 0^+} (f|\alpha)(x+iy)=\sumn b_n e(nx/h).\end{equation}
The distributions $\tau$ and $\tau|\alpha$ inherit the relation (\ref{intslash}), hence
\begin{equation}\label{intslashdist}
    (\tau|\alpha)(x)=j(\alpha,x)^{-k}\tau(\alpha x).
\end{equation}
\begin{lem}\label{tauintinforder}
The automorphic distribution $\tau$ attached to $f\in\Gamma_0(N)$ vanishes to infinite order at every rational number $\frac{p}{q}$. In particular, $\tau$ can be expressed, locally near $\frac{p}{q}$, as the sum
\begin{equation}\label{tausumlevelN}
\tau(x)=\sum_{j=0}^{m+k+2} \frac{d^j}{dx^j}\left(x-\frac{p}{q}\right)^{m+j} h_j(x),\end{equation} where $h_j$ are continuous functions on a neighborhood of $\frac{p}{q}$.
\end{lem}
\begin{proof}
Fix $\frac{p}{q}\in\Q$. Let $\gamma\in \SL(2,\Z)$ be a matrix such that $\gamma (i\infty)=\frac{p}{q}$. Then $\gamma=\begin{pmatrix}
p & r\\ 
q & \tilde{p}
\end{pmatrix}$, where $\tilde{p}$ is a multiplicative inverse of $p$ modulo $q$ and $r=(p\tilde{p}-1)/q$. By (\ref{intslashdist}), we have
\begin{equation}
    (\tau|\gamma)(\gamma^{-1}x)=j(\gamma,\gamma^{-1}x)^{-k}\tau(x)=j(\gamma^{-1},x)^k \tau(x).
\end{equation}
It follows that
\begin{equation}
    \tau(x)=j(\gamma^{-1},x)^{-k}(\tau|\gamma)(\gamma^{-1}x),
\end{equation}
hence we can write
\begin{equation}\label{tauintchain}
\tau(x)=(-qx+p)^{-k}(\tau|\gamma)\left(\frac{\tilde{p}x-r}{-qx+p}\right)=(-qx+p)^{-k}\left((-qx+p)\frac{d}{dx}\right)^l (\tau|\gamma)^{(-l)}\left(\frac{\tilde{p}x-r}{-qx+p}\right)
\end{equation}
for any $l>0$. Following the argument in Lemma~\ref{tauinforderlem}, it follows from (\ref{tauintchain}) that for any $m>0$ we can write $\tau$ as the sum
\begin{equation}\tau(x)=\sum_{j=0}^{m+k+2} \frac{d^j}{dx^j}\left(x-\frac{p}{q}\right)^{m+j} h_j(x),\end{equation}
where $h_j$ are continuous functions on a neighborhood of $p/q$.
\end{proof}
For a rational number $\frac{p}{q}\in\Q$, we fix the denominator $q$ to be positive and define \begin{equation}\label{taupqint}
    \tau_{p/q}(x)=\tau\left(\frac{x}{q}+\frac{p}{q}\right).
\end{equation} We derive the functional equation of $L_{p/q}(s)$ using the Mellin transform of $\tau_{p/q}$. 
\begin{lem}\label{intfunceqlem}
Let $\frac{p}{q}$ with $q>0$ be a rational number which is $\Gamma_0(N)$-equivalent to $i\infty$. Then the even Mellin transform of $\tau_{p/q}$ is an entire function, and satisfies the functional equations
\begin{equation}\label{intfunceq1}
M_0\tau_{p/q}(s)=q^sL_{p/q}\left(s-\frac{k-1}{2}\right)G_0(s),\end{equation}
and
\begin{equation}\label{intfunceq2}
    M_0\tau_{p/q}(s)=M_0\tau_{-\tilde{p}/q}(k-s).
\end{equation}
\end{lem}
\begin{proof}
    Taking the even Mellin transform of $\tau_{p/q}$ gives
\begin{equation}
M_0\tau_{p/q}(s)=\inR \tau\left(\frac{x}{q}+\frac{p}{q}\right)\abs{x}^{s-1}\dx =q^s\inR \tau\left(x+\frac{p}{q}\right)\abs{x}^{s-1}\dx.
\end{equation}
The distribution $\tau(x+p/q)$ is periodic, and can be written as a Fourier series
\begin{equation}
    \tau\left(x+\frac{p}{q}\right)=\sumn a_n n^{\frac{k-1}{2}} e^{2\pi i n \frac{p}{q}}e(nx).
\end{equation}
The functional equation (\ref{intfunceq1}) follows from Lemma~\ref{Mellinrightlem}. It follows from Lemma~\ref{tauintinforder} that $\tau\left(x+\frac{p}{q}\right)$ vanishes to infinite order at $0$, and therefore by Lemma~\ref{Mellinentirelem} the Mellin transform is entire. \par 
Let $\tilde{p}$ be a multiplicative inverse of $p$ mod $q$ and let $\gamma^{-1}=\begin{pmatrix}
\tilde{p} & -r \\ 
-q & p
\end{pmatrix}\in \Gamma_0(N)$. We have
\begin{equation}\label{taupqaut}
    \tau\left(\frac{x}{q}+\frac{p}{q}\right)=(-x)^{-k}\tau\left(\begin{pmatrix}
\tilde{p} & -r \\ 
-q & p
\end{pmatrix}\left(\frac{x}{q}+\frac{p}{q}\right)\right)=\abs{x}^{-k}\tau\left(-\frac{1}{qx}-\frac{\tilde{p}}{q}\right).
\end{equation}
The automorphic relation (\ref{taupqaut}) and the change of variable from $x$ to $-x^{-1}$ yield the identity
\begin{equation}\label{intfunceq3}
\begin{aligned}
    M_0\tau_{p/q}(s)&=\inR \tau\left(\frac{x}{q}+\frac{p}{q}\right)\abs{x}^{s-1}\dx=\inR \tau\left(-\frac{1}{qx}-\frac{\tilde{p}}{q}\right)\abs{x}^{-k-s-1}\dx \\
    &=\inR \tau\left(\frac{x}{q}-\frac{\tilde{p}}{q}\right)\abs{x}^{k-s-1}\dx=M_0\tau_{-\tilde{p}/q}(k-s).
\end{aligned}\end{equation}
The change of variable $x$ to $-1/x$ is legitimate, since $\tau$ vanishes to infinite order both at $-\frac{\tilde{p}}{q}$ and at $\infty$.
\end{proof}
\begin{cor}\label{intrealvalued}
Let $f\in S_k\left(\Gamma_0(N)\right)$ be a cusp form such that the Fourier coefficients of $f$ are all real-valued. If $\frac{p}{q}$ is a rational number which is $\Gamma_0(N)$-equivalent to $i\infty$ and satisfies $p^2\equiv 1\pmod q$, then $M_0\tau_{p/q}(k/2+it)$ is real-valued for all $t\in\R$, and the function
\[Z_f(t)=i^{-k/2} \left(\frac{2\pi}{q}\right)^{k/2}\frac{M_0\tau_{p/q}\left(\frac{k}{2}+it\right)}{2\cos\left(\frac{\pi}{2} (k/2+it)\right)\abs{\Gamma(k/2+it)}}\]
is a real-valued function on the real line such that ${Z_f(t)}=\abs{L_{p/q}\left(\frac{1}{2}+it\right)}$.
\end{cor}
\begin{proof}
Since $p$ is its own multiplicative inverse modulo $q$, we have $\tilde{p}=p$. The functional equation (\ref{intfunceq2}) with $s=\frac{k}{2}+it$, $t\in\R$ gives
\begin{equation}\label{funceqcrit}
    M_0\tau_{p/q}(k/2+it)= M_0\tau_{-\tilde{p}/q}(k/2-it)=M_0\tau_{-p/q}(k/2-it).
\end{equation}
From (\ref{intfunceq1}) and the definition (\ref{distgamma}), we have
\begin{equation}\label{funceqcritright}
    M_0\tau_{-p/q}(k/2-it)=2\left(\frac{2\pi}{q}\right)^{-\frac{k}{2}+it}L_{-p/q}(k/2-it)\Gamma(k/2-it)\cos\left(\frac{\pi}{2}(k/2-it)\right)
\end{equation}
and
\begin{equation}\label{funceqcritleft}
   M_0\tau_{p/q}(k/2+it)=2\left(\frac{2\pi}{q}\right)^{-\frac{k}{2}-it}L_{p/q}(k/2+it)\Gamma(k/2+it)\cos\left(\frac{\pi}{2}(k/2+it)\right). 
\end{equation}
Since all of the Fourier coefficients of $f$ are real numbers, they are complex conjugate of each other for $t\in\R$. Hence $M_0\tau_{p/q}(k/2+it)$ is real-valued for $t\in\R$. The identity
\[(-1)^n\cos(z)=\cos(n\pi-z),\quad n\in\Z\]
with $n=k/2$ and $z=\frac{\pi}{2}(k/2+it)$ implies that 
\[i^{k/2}\cos\left(\frac{\pi}{2}(k/2+it)\right)=i^{-k/2}\cos\left(\frac{\pi}{2}(k/2-it)\right),\]
hence $i^{k/2}\cos\left(\frac{\pi}{2}(k/2+it)\right)$ is real-valued for $t\in\R$. Therefore the function
\[Z_f(t)=i^{-k/2} \left(\frac{2\pi}{q}\right)^{k/2}\frac{M_0\tau_{p/q}\left(\frac{k}{2}+it\right)}{2\cos\left(\frac{\pi}{2} (k/2+it)\right)\abs{\Gamma(k/2+it)}}\]
is real-valued for all $t\in\R$, and it is the direct consequence of the functional equation (\ref{funceqcritleft}) that $\abs{Z_f(t)}=\abs{L_{p/q}\left(\frac{1}{2}+it\right)}$.
\end{proof}

\section{Modular forms of half-integral weight}
In this section we define the automorphic distributions attached to cusp forms of half-integral weight. We define $\sqrt z=z^{1/2}$ so that $-\pi/2<\arg (z^{1/2})\leq \pi/2$. For a positive integer $k$ and $N$ a positive multiple of $4$, we let $S_{k+\frac{1}{2}}\left(\Gamma_0(N)\right)$ denote the space of cusp forms of weight $k+\frac{1}{2}$ for the congruence subgroup $\Gamma_0(N)$. We refer to Shimura \cite{shimurahalfcusp} for classical facts on modular forms of half-integral weight. A cusp form $g(z)\in S_{k+\frac{1}{2}}\left(\Gamma_0(N)\right)$ satisfies
\begin{equation}\label{halfautomorphy}
    g\left(\frac{az+b}{cz+d}\right)=\left(\frac{c}{d}\right)^{2k+1}\epsilon_d^{-1-2k}(cz+d)^{k+\frac{1}{2}}g(z) 
\end{equation}
for any $\begin{pmatrix}
a &b \\ 
 c& d
\end{pmatrix}\in \Gamma_0(N)$, where $\epsilon_d$ is $1$ or $i$ for odd $d$ according to whether $d\equiv 1\pmod 4$ or $d\equiv 3\pmod 4$, respectively, and $\left(\frac{c}{d}\right)$ is Shimura's extension of the Jacobi symbol.
We define the automorphic distribution $\tau=\tau_g$ attached to $g$ to be $\tau(x)=\lim_{y\to 0^+}g(x+iy)$.
The automorphy relation (\ref{halfautomorphy}) with $c=0$ gives $g(z+1)=g(z)$, hence $g(z)$ has a Fourier expansion of the form
\begin{equation}\label{gfourier}
    g(z)=\sumn b_n n^{\frac{k}{2}-\frac{1}{4}}e(nz).
\end{equation}
It follows that $\tau(x+1)=\tau(x)$, and $\tau$ can also be written as a Fourier series
\begin{equation}\label{taugfourier}
\tau(x)=\sumn b_n n^{\frac{k}{2}-\frac{1}{4}}e(nx).
\end{equation}
For $x,y\in \R_{\neq 0}$, let $(a,b)_H$ denote the Hilbert symbol for $\R$ given by the formula
\begin{equation}\label{hilbertsym}
    (x,y)_H=\begin{cases}
        -1 &\text{if } x<0 \text{ and }y<0 \\
        1 &\text{otherwise}.
    \end{cases}
\end{equation}
Observe that for $c\neq 0$ and $cx+d\neq 0$ we have
\begin{equation}\label{sqrtlim}
\lim_{y\to 0^+} \sqrt{c(x+iy)+d}=
(c,cx+d)_H\sgn(cx+d)\sqrt{\abs{cx+d}}.
\end{equation}
Thus the distribution $\tau$ satisfies the automorphy relation
\begin{equation}\label{taugautomorphy}
     \tau\left(\frac{ax+b}{cx+d}\right)=
     \left(\frac{c}{d}\right)^{2k+1}\epsilon_d^{-1-2k}(c,cx+d)_H\sgn(cx+d)^{k+\frac{1}{2}}\abs{cx+d}^{k+\frac{1}{2}}\tau(x),
\end{equation}
for $\begin{pmatrix}
a & b\\ 
c & d
\end{pmatrix}\in\Gamma_0(N)$ with $c\neq 0$. \par
The slash operators on the space of half-integral cusp forms are defined as follows. Let $\mathcal{G}$ denote the set of all couples $(\alpha,\phi)$ formed by an element $\alpha=\begin{pmatrix}
a & b\\ 
c & d
\end{pmatrix}\in \GL^{+}(2,\R)$ and a holomorphic function $\phi(z)$ on the upper half plane such that $\phi(z)^2=t\cdot \det(\alpha)^{-1/2}(cz+d)$ with $\abs{t}=1$. Then $\mathcal G$ forms a group with the multiplication law 
\begin{equation}
\left(\alpha,\phi(z)\right) \left(\beta,\psi(z)\right)=\left(\alpha\beta,\phi(\beta z)\psi(z)\right).\end{equation} 
For $\xi=\left(\alpha,\phi(z)\right)\in\mathcal G$, we define the slash operator as
\begin{equation}\label{halfslash}
(g|\xi) (z)=g(\alpha z) \phi(z)^{-2k-1}.\end{equation} 
As in the integral weight case, we shall show that $\tau$ vanishes to infinite order at every cusp using the slash operator. Fix a rational number $\frac{p}{q}\in\Q$ with $q>0$. Let $\gamma=\begin{pmatrix}
p & r\\ 
q & \tilde{p}
\end{pmatrix}\in\SL(2,\Z)$ and $\phi(z)=\sqrt{qz+\tilde{p}}$. Then $\rho=\left(\gamma,\phi(z)\right)$ is an element of $\mathcal G$. Since $g$ is a cusp form, the function $g|\rho$ can be written as a Fourier series
\begin{equation}
    (g|\rho)(z)=\sum_{n\geq 0} c_n e\left((n+r)z/h\right)
\end{equation}
for some $0\leq r<1$, and $c_0=0$ in case $r=0$. The value of $r$ depends only on the $\Gamma_0(N)$-equivalence class of $\frac{p}{q}$ (see \cite{shimurahalfcusp} for details). Let $\tau|\rho$ be the distribution which corresponds to $g|\rho$:
\begin{equation}
    (\tau|\rho)(x)=\lim_{y\to 0^+}(g|\rho)(x+iy)=\sum_{n\geq 0} c_n e\left((n+r)x/h\right).
\end{equation}
As in (\ref{taujth}), the distribution $\tau|\rho$ can be written as the $j$-th derivative of a continuous periodic function for every sufficiently large integer j,
\begin{equation}\label{tauslashjth}
    (\tau|\rho)^{(-j)}(x)=\sum_{n\geq 0} \left(2\pi i (n+r)/h\right)^{-j}c_n e\left((n+r)x/h\right).
\end{equation}
From (\ref{sqrtlim}) and (\ref{halfslash}), we see that $\tau$ and $\tau|\rho$ satisfy
\begin{equation}\label{halfslashdist}
    (\tau|\rho)(x)=\tau(\gamma x)(q,qx+\tilde{p})_H^{-2k-1}\sgn(qx+\tilde{p})^{-2k-1}\sqrt{\abs{qx+\tilde{p}}}^{-2k-1}.
\end{equation}
Observe that
\[q(\gamma^{-1} x)+\tilde{p}=q\left(\frac{px-r}{-qx+p}\right)+\tilde{p}=\frac{1}{-qx+p}.\]
Hence
\begin{equation}\label{halfslashdist2}
\begin{aligned}
   (\tau|\rho)(\gamma^{-1}x)&=\tau( x)(q,(-qx+p)^{-1})_H^{-2k-1}\sgn((-qx+p)^{-1})^{-2k-1}\sqrt{\abs{-qx+p}}^{2k+1}\\ 
   &=\tau( x)(q,-qx+p)_H\sgn(-qx+p)\sqrt{\abs{-qx+p}}^{2k+1}.  
\end{aligned}
\end{equation}
\begin{lem}\label{tauhalfinforder}
The automorphic distribution $\tau$ attached to $g\in\Gamma_0(N)$ vanishes to infinite order at every rational number $\frac{p}{q}$. In particular, $\tau$ can be expressed, locally near $\frac{p}{q}$, as the sum
\begin{equation}\label{tausumlevelNhalf}
\tau(x)=\sum_{j=0}^{m+k+3} \frac{d^j}{dx^j}\left(x-\frac{p}{q}\right)^{m+j} h_j(x),\end{equation} where $h_j$ are continuous functions on a neighborhood of $\frac{p}{q}$.
\end{lem}
\begin{proof}
Fix a rational number $\frac{p}{q}\in\Q$ with $q>0$. Let $\rho=\left(\gamma,\phi(z)\right)\in \mathcal G$ with $\gamma=\begin{pmatrix}
p & r\\ 
q & \tilde{p}
\end{pmatrix}$ and $\phi(z)=\sqrt{qz+\tilde{p}}$. Applying the chain rule to the equation (\ref{halfslashdist2}), we obtain
\begin{equation}\label{halfslashdist3}
    \begin{aligned}
   \tau( x)&=(q,-qx+p)_H\sgn(-qx+p)\abs{-qx+p}^{-k-\frac{1}{2}}(\tau|\rho)(\gamma^{-1}x) \\
   &=(q,-qx+p)_H\sgn(-qx+p)\left((-qx+p)\frac{d}{dx}\right)^l \left(\tau|\rho\right)^{(-l)}\left(\frac{\tilde px-r}{-qx+p}\right).
    \end{aligned}
\end{equation}
Following the argument in Lemma~\ref{tauinforderlem} with $l=m+k+2$, we can remove the pole of $\abs{-qx+p}^{-k-\frac{1}{2}}$ at $\frac{p}{q}$ and obtain additional vanishing of order $m$. We conclude that $\tau$ can be written as
\begin{equation}
    \tau(x)=\sum_{j=0}^{m+k+3} \frac{d^j}{dx^j}\left(x-\frac{p}{q}\right)^{m+j}h_j(x),
\end{equation}
locally near $\frac{p}{q}$, where $h_j$ are continuous functions on a neighborhood of $\frac{p}{q}$.
\end{proof}
For a rational number $\frac{p}{q}$ with $q>0$, we define 
\begin{equation}\label{taupq}
    \tau_{p/q}(x)=\tau\left(\frac{x}{q}+\frac{p}{q}\right).
\end{equation}
We derive the functional equation of $L_{p/q}(s)$ using the Mellin transform of $\tau_{p/q}$.
\begin{lem}\label{halffunceqlem2}
Let $\frac{p}{q}$ with $q>0$ be a rational number which is $\Gamma_0(N)$-equivalent to $i\infty$. Then the even Mellin transform of $\tau_{p/q}$ is an entire function, and satisfies the functional equations
\begin{equation}\label{halffunceq1}
M_0\tau_{p/q}(s)=q^s L_{p/q}\left(s-\frac{k}{2}+\frac{1}{4}\right)G_0(s),\end{equation}
and
\begin{equation}\label{halffunceq2}
     \frac{M_0\tau_{p/q}(s)}{2\cos(\pi s/2)}=i^{k+\frac{1}{2}}\left(\frac{-q}{p}\right)^{-1-2k}\epsilon_p^{2k+1}\left(\frac{2\pi}{q}\right)^{-k-\frac{1}{2}+s}\Gamma\left(k+\frac{1}{2}-s\right)L_{-\tilde{p}/q}\left(\frac{k}{2}+\frac{3}{4}-s\right).
\end{equation}
\end{lem}
\begin{proof}
As in Lemma~\ref{intfunceq1}, the functional equation (\ref{halffunceq1}) follows from Lemma~\ref{Mellinrightlem}, and it follows from Lemma~\ref{Mellinentirelem} and Lemma~\ref{tauhalfinforder} that the Mellin transform is entire. Let $\tilde{p}$ be a multiplicative inverse of $p$ mod $q$ and let $\gamma^{-1}=\begin{pmatrix}
\tilde{p} & -r \\ 
-q & p
\end{pmatrix}\in \Gamma_0(N)$. From the automorphy relation (\ref{taugautomorphy}) with $\gamma^{-1}$, we have
\begin{equation}\label{taupqhalfaut}
    \tau\left(\frac{x}{q}+\frac{p}{q}\right)= \left(\frac{-q}{p}\right)^{-1-2k}\epsilon_p^{2k+1}(-q,-x)_H\sgn(-x)^{-k-\frac{1}{2}}\abs{x}^{-k-\frac{1}{2}}\tau\left(-\frac{1}{qx}-\frac{\tilde{p}}{q}\right).
\end{equation}
By (\ref{taupqhalfaut}), we have
\begin{equation}\label{halffunceq3}
\begin{aligned}
    M_0\tau_{p/q}(s) & = \inR \tau\left(\frac{x}{q}+\frac{p}{q}\right)\abs{x}^{s-1}\dx \\ 
    &=\left(\frac{-q}{p}\right)^{-1-2k}\epsilon_p^{2k+1}\inR (-q,-x)_H\sgn(-x)^{-k-\frac{1}{2}}\tau\left(-\frac{1}{qx}-\frac{\tilde{p}}{q}\right)\abs{x}^{-k+s-\frac{3}{2}}\dx \\
    &=\left(\frac{-q}{p}\right)^{-1-2k}\epsilon_p^{2k+1}\inR (-q,x^{-1})_H\sgn(x^{-1})^{-k-\frac{1}{2}}\tau\left(\frac{x}{q}-\frac{\tilde{p}}{q}\right)\abs{x}^{k-s-\frac{1}{2}}\dx \\
    &=\left(\frac{-q}{p}\right)^{-1-2k}\epsilon_p^{2k+1}\inR \sgn(x)^{-k+\frac{1}{2}}\tau\left(\frac{x}{q}-\frac{\tilde{p}}{q}\right)\abs{x}^{k-s-\frac{1}{2}}\dx.
\end{aligned}\end{equation}
The last equality uses the fact that $-q<0$, hence $(-q,x^{-1})_H=\sgn(x^{-1})=\sgn(x)$. The integral on the right evaluates as
\begin{equation}\label{halffunceq4}
    \begin{aligned}
&\inR \sgn (x)^{-k+\frac{1}{2}}\tau\left(\frac{x}{q}-\frac{\tilde{p}}{q}\right) \abs{x}^{k-s-\frac{1}{2}}\dx \\
&\quad=q^{k+\frac{1}{2}-s}\inR \sgn(x)^{-k+\frac{1}{2}}\tau\left(x-\frac{\tilde{p}}{q}\right) \abs{x}^{k-s-\frac{1}{2}}\dx \\
&\quad=q^{k+\frac{1}{2}-s}\sumn a_n e(-n \tilde{p}/q)n^{\frac{k}{2}-\frac{1}{4}}\inR \sgn(x)^{-k+\frac{1}{2}}e(nx) \abs{x}^{k-s-\frac{1}{2}}\dx \\
&\quad = q^{k+\frac{1}{2}-s}L_{-\tilde{p}/q}\left(\frac{k}{2}+\frac{3}{4}-s\right)\inR\sgn(x)^{-k+\frac{1}{2}} e(x) \abs{x}^{k-s-\frac{1}{2}}\dx. \end{aligned}\end{equation}
As in Lemma~\ref{Mellinrightlem} and Lemma~\ref{Mellinentirelem}, the fact that $\tau$ vanishes to infinite order at $-\frac{\tilde p}{q}$ and at $\infty$ legitimizes the interchange of the order of integration and summation. \par 
The functional equation (\ref{halffunceq2}) follows once we show that
\begin{equation}\label{halfgamma}
\inR\sgn(x)^{-k+\frac{1}{2}} e(x) \abs{x}^{k-s-\frac{1}{2}}\dx=2i^{k+\frac{1}{2}}(2\pi)^{-k-\frac{1}{2}+s}\cos\left(\frac{\pi s}{2}\right)\Gamma\left(k+\frac{1}{2}-s\right).
\end{equation}
For $0<\real \nu<1$ the integral
\[
    \int_{-m_1}^{m_2} e(x)\abs{x}^{\nu-1}\sgn(x)^{-k+\frac{1}{2}}\dx=\int_{0}^{m_2}e(x)\abs{x}^{\nu-1}\dx+(-1)^{-k+\frac{1}{2}}\int_{0}^{m_1}e(-x)\abs{x}^{\nu-1}\dx
\]
converges to $(2\pi)^{-\nu}\Gamma(\nu)\left(i^\nu+i^{2k+1-\nu}\right)$
as $m_1,m_2\to\infty$. Invoking the uniqueness of meromorphic continuation with $\nu=k+\frac{1}{2}-s$ deduces (\ref{halfgamma}).
\end{proof}
\begin{cor}\label{halfrealvalued}
Let $g\in S_{k+\frac{1}{2}}\left(\Gamma_0(N)\right)$ be a cusp form such that the Fourier coefficients of $g$ are all real-valued. If $\frac{p}{q}$ is a rational number which is $\Gamma_0(N)$-equivalent to $i\infty$ and satisfies $p^2\equiv 1\pmod q$, then the function
\[H_g(t)=\beta_{p/q}^{\frac{1}{2}}i^{-\frac{k}{2}-\frac{1}{4}}\frac{M_0\tau_{p/q}\left(\frac{k}{2}+\frac{1}{4}+it\right)}{2\cos\left(\frac{\pi}{2} (k/2+1/4+it)\right)},\]
where $\beta_{p/q}=\left(\frac{-q}{p}\right)^{2k+1}\epsilon_p^{-1-2k}$,
is real-valued for all $t\in\R$. The function
\[Z_g(t)=\beta_{p/q}^{\frac{1}{2}}i^{-\frac{k}{2}-\frac{1}{4}} \left(\frac{2\pi}{q}\right)^{\frac{k}{2}+\frac{1}{4}}\frac{M_0\tau_{p/q}\left(\frac{k}{2}+\frac{1}{4}+it\right)}{2\cos\left(\frac{\pi}{2} (k/2+1/4+it)\right)\abs{\Gamma(k/2+1/4+it)}}\]
is a real-valued function on the real line such that $\abs{Z_g(t)}=\abs{L_{p/q}\left(\frac{1}{2}+it\right)}$.
\end{cor}
\begin{proof}
By the functional equations (\ref{halffunceq1}) and (\ref{halffunceq2}), we have
\begin{equation}\label{Zgt1}
\begin{aligned}
   \frac{M_0\tau_{p/q}\left(\frac{k}{2}+\frac{1}{4}+it\right)}{2\cos\left(\frac{\pi}{2} (k/2+1/4+it)\right)}&= \left(\frac{2\pi}{q}\right)^{-\frac{k}{2}-\frac{1}{4}-it}\Gamma\left(\frac{k}{2}+\frac{1}{4}+it\right)L_{p/q}\left(\frac{1}{2}+it\right) \\
   &=\beta_{p/q}^{-1}i^{k+\frac{1}{2}}\left(\frac{2\pi}{q}\right)^{-\frac{k}{2}-\frac{1}{4}+it}\Gamma\left(\frac{k}{2}+\frac{1}{4}-it\right)L_{-\tilde{p}/q}\left(\frac{1}{2}-it\right).
\end{aligned} 
\end{equation}
If $p=\tilde{p}$ then the function $H_g(t)$ equal to its own complex conjugate. It follows that $Z_g(t)$ is a real-valued function such that $\abs{Z_g(t)}=\abs{L_{p/q}\left(\frac{1}{2}+it\right)}$.
\end{proof}

\section{Preliminary lemmas}
In this section we construct a parameterized family of test functions $u_T$ with the desired properties described in Section~\ref{outline}, and state some lemmas that will be used in the proof of Theorem~\ref{integralthm} and Theorem~\ref{halfintegralthm}. \par
Let $\varphi(x)$ be an even bump function supported on $[-1,1]$ satisfying $0\leq \varphi(x)\leq 1$ and $\varphi(0)=1$. Let 
\begin{equation}\label{defpsi}
    \psi(x)=\varphi(\log x).
\end{equation} 
Then $\psi(x)=\psi(1/x)$ and $\psi(x)$ is supported on $[1/e,e]$.
Next, let $\lambda(x)$ be the convolution of $\psi(x)$ with itself.
\begin{equation}\label{deflambda}
\lambda(x)=(\psi*\psi)(x)=\inP \psi(y)\psi \left(\frac{x}{y}\right)\frac{dy}{y}.
\end{equation}
Then $\lambda(x)$ is a smooth function supported on $[1/e^2,e^2]$, and satisfies $\lambda(x)=\lambda(1/x)$. For $T>0$, define the parameterized family of functions $\lambda_T(x)$ as 
\begin{equation}\label{deflambdaT}
\lambda_T(x)=Tx^{-2iT^{\frac{3}{2}}}\lambda(x^T).\end{equation}
The functions $\lambda_T(x)$ are supported on $[e^{-2/T},e^{2/T}].$
For $t\in\R$ and $T>0$, we let $u(t)$ and $u_T(t)$ be the Mellin transforms
\begin{equation}\label{defu}
u(t)=M\lambda(it) \end{equation} and
\begin{equation}\label{defuT}
        u_T(t)=M\lambda_T(it).
\end{equation}
Observe that
\begin{equation}\label{uTtranslate}
\begin{aligned}
    u_T(t)&=\inP \lambda_T(x)x^{it}\,\frac{dx}{x}=\inP T\lambda(x^T)x^{it-2iT^{\frac{3}{2}}}\,\frac{dx}{x}=\inP \lambda(x) x^{\frac{it-2iT^{\frac{3}{2}}}{T}}
    \,\frac{dx}{x} \\
    &=
    u\left(\frac{t-2T^{\frac{3}{2}}}{T}\right).
    \end{aligned}
\end{equation}
\begin{lem}
The function $u_T(t)$ is real-valued, and satisfies the inequality $0\leq u_T(t)\leq 4$ for all $t\in\R$ and $T>0$.
\end{lem}
\begin{proof}
By \ref{uTtranslate}, it suffices to show that $u(t)$ is real-valued and satisfies the inequality $0\leq u(t)\leq 4$ for all $t\in\R$. Since $\lambda(x)$ is the convolution of $\psi(x)$, we have \begin{equation}\label{utsq}
    u(t)=M\lambda(it)=\left(M\psi(it)\right)^2.
\end{equation} 
Since $\psi(x)=\psi(x^{-1})$, we have
\begin{equation}\label{Mpsiit}
M\psi(it)=\inP \psi(x)x^{it}\xdx=\int_0^{1} \psi(x)(x^{it}+x^{-it})\xdx=2\int_0^{1}\psi(x)\cos(t\log x)\xdx.\end{equation}
Hence $M\psi(it)$ is real-valued for all $t\in\R$, consequently $u(t)$ is real-valued and nonnegative for all $t\in\R$. The boundedness also follows immediately from (\ref{Mpsiit}), since $\psi(x)$ is supported away from $0$.
\end{proof}
\begin{lem}\label{uTboundlem}
For any $N>0$, the function $u_T(t)$ satisfies the bound \begin{equation}\label{uTbound}
u_T(t+2T^{3/2})\ll_N T^N\abs{t}^{-N}.\end{equation}
\end{lem}
\begin{proof}
Since $u_T(t+2T^{3/2})=u(t/T)$, it suffices to show that $u(t)\ll_N \abs{t}^{-N}$ as $\abs{t}\to \infty$. Since $u(t)=M\lambda(it)$ and $\lambda$ is a Schwartz function, the asymptotic follows from the well-known fact that the Mellin transform of a Schwartz function decays rapidly on vertical lines.
\end{proof}
The next two lemmas establish bounds of integrals of the product of $u_T(x)$ and polynomially bounded functions.
\begin{lem}\label{absfuTlem}
Let $f(x)$ be a continuous real-valued function on the real line that satisfies the bound $f(x)=O(\abs{x}^\alpha)$ as $\abs{x}\to\infty$ for some $\alpha>0$. Then the integral
\begin{equation}\int_{\abs{x-2T^{3/2}}>T^{3/2}} \abs{f(x)}u_T(x)\dx.\end{equation}
decays rapidly as $T\to\infty$.
\end{lem}
\begin{proof}
For large $T>0$ and large $N>0$, we have the bound
\begin{equation}\begin{aligned}
    \int_{\abs{x-2T^{3/2}}>T^{3/2}} \abs{f(x)}u_T(x)\dx &=\int_{\abs{x}>T^{3/2}} \abs{f(x+2T^{3/2})}u_T(x+2T^{3/2})\dx \\
    &\ll_N \int_{\abs{x}>T^{3/2}} (\abs{x}+2T^{3/2})^\alpha T^N \abs{x}^{-N} \dx \\
    &\ll_N 3^\alpha T^N \int_{\abs{x}>T^{3/2}} \abs{x}^{a-N}\dx \\
    &\ll_{\alpha,N} T^{-\frac{N}{2}+\frac{3a}{2}+\frac{3}{2}}.
\end{aligned}\end{equation}
The first inequality follows from Lemma~\ref{uTboundlem}. Therefore the integral decays rapidly as $T\to\infty$.
\end{proof}
\begin{lem}\label{fuTdifflem}
Let $f(x)$ be as defined in Lemma~\ref{absfuTlem}. If $f(x)$ has only finitely many zeros on the real line, then the difference
\begin{equation}\inR \abs{f(x)}u_T(x)\dx-\abs{\inR f(x)u_T(x)\dx}.\end{equation}
decays rapidly in $T$.
\end{lem}
\begin{proof}
If $f(x)$ has only finitely many zeros, then we can choose a large $T>0$ such that $f(x)$ does vanish in the interval $[T^{3/2},3T^{3/2}]$. For such $T>0$, we have
\begin{equation}\int_{T^{3/2}}^{3T^{3/2}}\abs{f(x)}u_T(x)\dx=\abs{\int_{T^{3/2}}^{3T^{3/2}}f(x)u_T(x)\dx}.\end{equation}
Thus for any large $N>0$ we have
\begin{equation}\begin{aligned}
    &\inR \abs{F(x)}u_T(x)\dx-\abs{\inR F(x)u_T(x)\dx} \\
    &\quad = \int_{\abs{x-2T^{3/2}}>T^{3/2}} \abs{F(x)}u_T(x)\dx+\abs{\int_{T^{3/2}}^{3T^{3/2}}F(x)u_T(x)\dx}-\abs{\inR F(x)u_T(x)\dx} \\
    &\quad \leq \int_{\abs{x-2T^{3/2}}>T^{3/2}} \abs{F(x)}u_T(x)\dx+\abs{\int_{\abs{x-2T^{3/2}}>T^{3/2}} F(x)u_T(x)\dx} \\
    &\quad \leq 2\int_{\abs{x-2T^{3/2}}>T^{3/2}} \abs{F(x)}u_T(x)\dx.
\end{aligned}\end{equation}
By Lemma~\ref{absfuTlem}, the last integral decays rapidly as $T\to\infty$.
\end{proof}

Let $\nu\in\frac{1}{2}\Z$ and let $f(z)\in S_\nu\left(\Gamma_0(N)\right)$ be a cusp form of integral or half-integral weight. Let $\tau(x)$ be the automorphic distribution attached to $f(z)$, and let $\tau_{p/q}(x)=\tau(x/q+p/q)$ for $\frac{p}{q}\in\Q$.

\begin{lem}\label{stirlinglem}
For $\sigma\in\R$, the function $G_0(s)=2(2\pi)^{-s}\Gamma(s)\cos(\pi s/2)$ satisfies the asymptotic formula
\begin{equation}\label{stirling}
    \abs{G_0(\sigma+it)}\sim\left(\frac{\abs{t}}{2\pi}\right)^{\sigma-1/2},\quad \abs{t}\to\infty.
\end{equation}
In particular, for $\sigma\geq 1/2$ we have
\begin{equation}\label{GObounded}
    \frac{1}{G_0(\sigma+it)}=O(1), \quad \abs{t}\geq 1.
\end{equation}
\end{lem}
\begin{proof}
The asymptotic formula is a direct consequence of the Stirling's formula for the $\Gamma$-function.
\end{proof}

\begin{lem}
The additively twisted $L$-function $L_{p/q}(s)$ is polynomially bounded on the critical line $\real s=\frac{1}{2}$. That is, there is an $\alpha>0$ such that $L_{p/q}(1/2+it)=O(\abs{t}^\alpha)$
\end{lem}
\begin{proof}
This is well known. The additively twisted $L$-functions satisfy the functional equations (\ref{intfunceq1}), (\ref{intfunceq2}), (\ref{halffunceq1}) and (\ref{halffunceq2}), hence the bound follows from Stirling's formula for the $\Gamma$-function and the Phragmen-Lindel\"{o}f principle.
\end{proof}
The next two lemmas establish the bound of integrals of $L$-functions and certain distributions paired against $u_T(x)$.
\begin{lem}\label{LuTlem}
Let $A(s)$ be a function given by a series $A(s)=\sumn a_n n^{-s}$ with $a_1\neq 0$ and $a_n=O(n^{\alpha})$ for some $\alpha>0$. For sufficiently large $T>0$, we have
\begin{equation}
    B_T(s)=\inR A(s+it)u_T(t)\dt=2\pi a_1 T
\end{equation}
for all $s\in\mathbb{C}$.
\end{lem}
\begin{proof}
Take $s\in\mathbb{C}$ with a large real part, so that $A(s)$ converges absolutely. By the Mellin inversion formula, we have
\begin{equation}
    \begin{aligned}
    B_T(s) &= \inR \sumn a_n n^{-s-it}M\lambda_T(it)\,dt =\sumn a_n n^{-s} \inR n^{-it} M\lambda_T(it) \,dt\\
    &=2\pi \sumn a_n n^{-s} \lambda_T(n).
\end{aligned}
\end{equation}
Recall that the function $\lambda_T(x)$ is supported on $[e^{-2/T},e^{2/T}]$. It follows that if $T>2/\log 2$ then $\lambda_T(n)=0$ for all positive integers $n$ greater than $1$. It follows by the definition of $\lambda_T$ given in (\ref{deflambdaT}) that
\begin{equation} 
B_T(s)=2\pi a_1 \lambda_T(1)=2\pi a_1 \lambda(1) T=2\pi a_1 T.
\end{equation}
The lemma follows from the uniqueness of analytic continuation. 
\end{proof}
\begin{lem}\label{MuTboundlem}
Let $\sigma\in C^{-\infty}(\R)$ be a periodic distribution without constant term that vanishes to infinite order at $0$ and at $\pm 1$. Suppose that there exists an $N_0>0$ such that for every $m>0$ the distribution $\sigma$ can be written as a sum
\begin{equation}\label{sigmapresentation}
    \sigma(x)+\sigma(-x)=\sum_{j=0}^{N} \frac{d^j}{dx^j}\left((x-1)^{m+j}h_j(x)\right),
\end{equation}
for some $N\leq m+N_0$, where $I_1$ is a neighborhood of $1$ and $h_j$ are continuous functions on $I_1$. Then for any $k>0$ we have
\begin{equation}\label{Parsevalbound}
    \inR M_0\sigma(k+ix)u_T(x)\dx \ll_{k,N} T^{-N}.
\end{equation}
\end{lem}
\begin{proof}
Fix a large positive integer $m$. Then $\sigma$ has a presentation (\ref{sigmapresentation}) on a neighborhood $I_1$ of $1$, with $N\leq m+N_0$. Take a large $T>0$ such that $[e^{-\frac{2}{T}}, e^{\frac{2}{T}}]\subseteq I_1$, so the support of $\lambda_T(x)$ is contained in $I_1$. First, observe that since $\lambda(x)=\lambda(1/x)$ we have
\begin{equation}\label{MlambdaT}
    M\lambda_T(-s)=\inP Tx^{-2iT^{\frac{3}{2}}}\lambda(x^T) x^{-s-1}\dx=\inP Tx^{2iT^{\frac{3}{2}}}\lambda(x^T) x^{s-1}\dx.
\end{equation}
Recall that a periodic distribution without constant term vanishes to infinite order at $\infty$, so $\sigma$ vanishes to infinite order at $0$ and at $\infty$. Applying (\ref{MlambdaT}) and Corollary~\ref{Parsevalcor}, we obtain
\begin{equation}\label{Parsevalbound1}
\begin{aligned}
   &\inR M_0\sigma(k+ix)u_T(x)\dx = \inR M_0\sigma(k+ix)M\lambda_T(ix)\dx \\
   &\quad=\frac{1}{i}\inIm  M_0\sigma(k-s)M\lambda_T(-s)\ds=\frac{1}{i}\inIm  M_0\sigma(k-s)M\left(Tx^{2iT^{\frac{3}{2}}}\lambda(x^T)\right)(s)\ds \\
   &\quad=2\pi\inR \sigma(x)T\lambda(\abs{x}^T)\abs{x}^{2iT^{\frac{3}{2}}+k-1}\dx=2\pi\intres \left(\sigma(x)+\sigma(-x)\right)T\lambda(x^T)x^{2iT^{\frac{3}{2}}+k-1}\dx.
\end{aligned}
\end{equation}
To simplify notation, we write $g_T(x)=T\lambda(x^T)x^{2iT^{\frac{3}{2}}+k-1}$. Using integration by parts, we write
\begin{equation}\label{sigmagT}
    \begin{aligned}
    \intres\left(\sigma(x)+\sigma(-x)\right)g_T(x)\dx &=\intres \left(\sum_{j=0}^{N} \frac{d^j}{dx^j}\left((x-1)^{m+j}h_j(x)\right)\right)g_T(x)\dx \\
    &=\sum_{j=0}^{N}\intres\left(\frac{d^j}{dx^j}\left((x-1)^{m+j}h_j(x)\right)\right)g_T(x)\dx \\
    &=\sum_{j=0}^{N}(-1)^j\intres \left((x-1)^{m+j}h_j(x)\right)\frac{d^j}{dx^j}\left(g_T(x)\right)\dx.
    \end{aligned}
\end{equation}
Next, we estimate the size of the derivatives of $g_T(x)=Tx^{2iT^{\frac{3}{2}}+k-1}\lambda(x^T)$. Since $\lambda(x)$ is a smooth function of compact support, we have
\begin{equation}\frac{d^j}{dx^j}\lambda(x)\ll_j 1,\end{equation} 
and it implies
\begin{equation}
\frac{d^j}{dx^j}T\lambda(x^T) \ll_j T^{j+1}.\end{equation}
Applying the product rule, we obtain the bound
\begin{equation}
    \begin{aligned}
    \frac{d^j}{dx^j}\left(T\lambda(x^T)x^{2iT^{3/2}+k-1}\right)&=\sum_{i=0}^{j} \binom{j}{i} \frac{d^i}{dx^i}T\lambda(x^T)\frac{d^{(j-i)}}{dx^{(j-i)}}x^{2iT^{3/2}+k-1} \\
    &\ll_{k,j} \sum_{i=0}^{j} O(T^{i+1}T^{\frac{3}{2}(j-i)}) \\
    &\ll_{k,j} T^{\frac{3}{2}j+1}
    \end{aligned}
\end{equation}
on the neighborhood of $1$. From this bound and the fact that $N\leq m+N_0$, we deduce that each integral in the sum on the right of (\ref{sigmagT}) satisfies the bound
\begin{equation}
\begin{aligned}
    \intres \left((x-1)^{m+j}h_j(x)\right)\frac{d^j}{dx^j}\left(g_T(x)\right)\dx&\ll_{k,j} T^{-1}\cdot T^{-m-j} \cdot T^{\frac{3}{2}j+1} \\
    &\ll_{k,j}T^{-m+\frac{1}{2}j} \\ 
    &\ll_{k,j}T^{-\frac{m}{2}+\frac{N_0}{2}}.
\end{aligned}
\end{equation}
We conclude that
\begin{equation}
\begin{aligned}
   \inR M_0\sigma(k+ix)u_T(x)&=2\pi\intres \left(\sigma(x)+\sigma(-x)\right)g_T(x)\dx \\
   &=\sum_{j=0}^{N}(-1)^j\intres \left((x-1)^{m+j}h_j(x)\right)\frac{d^j}{dx^j}\left(g_T(x)\right)\dx \\
&\ll_{k,m}T^{-\frac{m}{2}+\frac{N_0}{2}}. 
\end{aligned}
\end{equation}
The lemma follows, since $\frac{m}{2}-\frac{N_0}{2}$ can be arbitrarily large.
\end{proof}

\section{Proof of the theorems}
In this section we use the lemmas we have developed to prove our main results.
\subsection{Proof of Theorem~\ref{integralthm}}
Let $f\in S_k\left(\Gamma_0(N)\right)$ be a cusp form such that all the Fourier coefficients $a_j$ of $f$ are real-valued and $a_1\neq 0$. Let $\tau=\tau_f$ be the automorphic distribution attached to $f$. Suppose that $\frac{p}{q}$ is a rational number that is $\Gamma_0(N)$-equivalent to $i\infty$ and satisfies $p^2\equiv 1 \pmod q$. Fix $q>0$, and let $\tau_{p/q}(x)=\tau(x/q+p/q)$. By Lemma~\ref{tauintinforder}, the distribution $\tau_{p/q}$ vanishes at $0$ and at $\pm 1$, and we can write $\tau_{p/q}(x)+\tau_{p/q}(-x)$ as
\begin{equation}
    \tau_{p/q}(x)+\tau_{p/q}(-x)=\sum_{j=0}^{m+k+2} \frac{d^j}{dx^j}\left(x-1
\right)^{m+j} h_j(x),\end{equation} where $h_j$ are continuous functions in a neighborhood of $1$. Let
\begin{equation}\label{Zftdef}
Z_f(t)=i^{-\frac{k}{2}} (2\pi)^{\frac{k}{2}}\frac{q^{it}M_0\tau_{p/q}\left(\frac{k}{2}+it\right)}{2\cos\left(\frac{\pi}{2} (k/2+it)\right)\abs{\Gamma(k/2+it)}}.
\end{equation}
By Corollary~\ref{intrealvalued}, the function $Z_f(t)$ is a real-valued function on the real line which satisfies $\abs{Z_f(t)}=\abs{L_{p/q}\left(\frac{1}{2}+it\right)}$. By Lemma~\ref{LuTlem}, we have
\begin{equation}
    \inR L_{p/q}\left(\frac{1}{2}+it\right)u_T(t)\dt=2\pi a_1 T,
\end{equation}
and this implies
\begin{equation}\label{ZftuTbound1}
    \inR \abs{Z_f(t)}u_T(t)\dt=\inR \abs{L_{p/q}\left(\frac{1}{2}+it\right)}u_T(t)\dt\gg T.
\end{equation}
On the other hand, by Lemma~\ref{MuTboundlem}, we have
\begin{equation}\label{intMuTbound}
   \inR M_0\tau_{p/q}\left(\frac{k}{2}+it\right) u_T(t)\dt\ll_{N} T^{-N}.
\end{equation}
Suppose, for a contradiction, that $L_{p/q}(1/2+it)$ has only finitely zeros. By Corollary~\ref{intrealvalued}, the function in the integral in (\ref{intMuTbound}) is real-valued for all $t\in\R$. By Lemma~\ref{fuTdifflem}, we have
\begin{equation}\label{intMuTbound1}
   \inR \abs{M_0\tau_{p/q}\left(\frac{k}{2}+it\right)} u_T(t)\dt\ll_{N} T^{-N}.
\end{equation}
Observe that
\begin{equation}\label{Zft1}
    \abs{Z_f(t)}=\frac{\abs{M_0\tau_{p/q}(k/2+it)}}{\abs{G_0(k/2+it)}}.
\end{equation}
It follows from Lemma~\ref{stirlinglem} that the bound (\ref{intMuTbound1}) implies
\begin{equation}\label{ZftuTbound2}
    \int_{\abs{t}\geq 1} \abs{Z_f(t)}u_T(t)\dt\ll_N T^{-N}.
\end{equation}
For $\abs{t}\leq 1$, we have
\begin{equation}\label{ZftuTbound3}
    \int_{\abs{t}\leq 1} \abs{Z_f(t)}u_T(t)\dt= \int_{\abs{t}\leq 1} \abs{L_{p/q}\left(\frac{1}{2}+it\right)}u_T(t)\dt=O(1).
\end{equation}
The bounds (\ref{ZftuTbound2}) and (\ref{ZftuTbound3}) together imply
\begin{equation}\label{ZftuTbound4}
    \inR \abs{Z_f(t)}u_T(t)\dt= O(1)
\end{equation}
as $T\to\infty$. This contradicts the bound (\ref{ZftuTbound1}).

\subsection{Proof of Theorem~\ref{halfintegralthm}}
Let $g\in S_{k+\frac{1}{2}}\left(\Gamma_0(N)\right)$ be a cusp form such that all the Fourier coefficients $b_j$ are real-valued and $b_1\neq 0$. Let $\tau=\tau_g$ be the automorphic distribution attached to $g$. Let $\frac{p}{q}\in\Q$ be a rational number that is $\Gamma_0(N)$-equivalent to $i\infty$ such that $p^2\equiv 1 \pmod q$. Fix $q>0$ and let $\tau_{p/q}(x)=\tau(x/q+p/q)$. By Lemma~\ref{tauhalfinforder}, the distribution $\tau_{p/q}$ vanishes at $0$ and at $\pm 1$, and we can write $\tau_{p/q}(x)+\tau_{p/q}(-x)$ as
\begin{equation}
    \tau_{p/q}(x)+\tau_{p/q}(-x)=\sum_{j=0}^{m+k+3} \frac{d^j}{dx^j}\left(x-1
\right)^{m+j} h_j(x),\end{equation} where $h_j$ are continuous functions in a neighborhood of $1$. Let
\begin{equation}\label{Hgtdef}
    H_g(t)=\beta_{p/q}^{\frac{1}{2}}i^{-\frac{k}{2}-\frac{1}{4}}\frac{M_0\tau_{p/q}\left(\frac{k}{2}+\frac{1}{4}+it\right)}{2\cos\left(\frac{\pi}{2} (k/2+1/4+it)\right)}
\end{equation}
and
\begin{equation}\label{Zgtdef}
    Z_g(t)=\beta_{p/q}^{\frac{1}{2}}i^{-\frac{k}{2}-\frac{1}{4}} \left(\frac{2\pi}{q}\right)^{\frac{k}{2}+\frac{1}{4}}\frac{M_0\tau_{p/q}\left(\frac{k}{2}+\frac{1}{4}+it\right)}{2\cos\left(\frac{\pi}{2} (k/2+1/4+it)\right)\abs{\Gamma(k/2+1/4+it)}},
\end{equation}
where $\beta_{p/q}=\left(\frac{-q}{p}\right)^{2k+1}\epsilon_p^{-1-2k}$.
By Lemma~\ref{LuTlem}, we have
\begin{equation}
    \inR L_{p/q}\left(\frac{1}{2}+it\right)u_T(t)\dt=2\pi b_1 T,
\end{equation}
and this implies
\begin{equation}\label{ZgtuTbound1}
    \inR \abs{Z_f(t)}u_T(t)\dt=\inR \abs{L_{p/q}\left(\frac{1}{2}+it\right)}u_T(t)\dt\gg T.
\end{equation}
On the other hand, by Lemma~\ref{MuTboundlem}, we have
\begin{equation}\label{gintMuTbound}
   \inR M_0\tau_{p/q}\left(\frac{k}{2}+it\right) u_T(t)\dt\ll_{N} T^{-N}.
\end{equation} 
From (\ref{Hgtdef}), we obtain the bound
\begin{equation}\label{HgtMuTbound}
\begin{aligned}
    &\inR 2\cos\left(\frac{\pi}{2}\left(k/2+1/4+it\right)\right) H_g(t) u_T(t) \dt \\
    &\quad=\inR \left(e\left(\frac{k}{8}+\frac{1}{16}\right)e^{-\frac{1}{2}\pi t}+e\left(-\frac{k}{8}-\frac{1}{16}\right)e^{\frac{1}{2}\pi t} \right) H_g(t)  u_T(t) \dt
    \ll_{N} T^{-N}.
\end{aligned}
\end{equation}
Observe that for any $\alpha,\beta,\theta_1,\theta_2\in\R$, the following inequality holds:
\begin{equation}\label{elemineq}
\abs{e^{i\theta_1}\alpha+e^{i\theta_2}\beta}\geq \abs{\sin(\theta_2-\theta_1)}\max\left(\abs{\alpha},\abs{\beta}\right)\geq \abs{\sin(\theta_2-\theta_1)}\frac{\abs{\alpha}+\abs{\beta}}{2}.
\end{equation}
Applying the inequality to (\ref{HgtMuTbound}), we obtain
\begin{equation}\label{HgtMuTbound2}
\begin{aligned}
    &\inR (e^{-\frac{1}{2}\pi t}+e^{\frac{1}{2}\pi t} ) H_g(t)  u_T(t) \dt
    \ll_{N} T^{-N}.
\end{aligned}
\end{equation}
Suppose, for a contradiction, that $L_{p/q}(1/2+it)$ has only finitely zeros. Observe that the function $(e^{-\frac{1}{2}\pi t}+e^{\frac{1}{2}\pi t}) H_g(t)$ is real-valued for all $t\in\R$. By Lemma~\ref{fuTdifflem}, we have
\begin{equation}\label{HgtMuTbound3}
\begin{aligned}
    &\inR (e^{-\frac{1}{2}\pi t}+e^{\frac{1}{2}\pi t} ) \abs{H_g(t)}  u_T(t) \dt
    \ll_{N} T^{-N}.
\end{aligned}
\end{equation}
Observe that
\begin{equation}\label{Zgtdef2}
    \abs{Z_g(t)}=\left(\frac{2\pi}{q}\right)^{\frac{k}{2}+\frac{1}{4}} \frac{(e^{-\frac{1}{2}\pi t}+e^{\frac{1}{2}\pi t})\abs{H_g(t)}}{(e^{-\frac{1}{2}\pi t}+e^{\frac{1}{2}\pi t})\abs{\Gamma(k/2+1/4+it)}}.
\end{equation}
The elementary bound
\begin{equation}
    \abs{\alpha}+\abs{\beta}\geq \abs{e^{i\theta_1}\alpha+e^{i\theta_2}\beta}
\end{equation}
implies
\begin{equation}
  e^{-\frac{1}{2}\pi t}+e^{\frac{1}{2}\pi t}\geq   2\cos\left(\frac{\pi}{2} (k/2+1/4+it)\right),
\end{equation}
hence it follows from Lemma~\ref{stirlinglem} that
\begin{equation}\label{stirlinghalf}
 \frac{1}{(e^{-\frac{1}{2}\pi t}+e^{\frac{1}{2}\pi t})\abs{\Gamma(k/2+1/4+it)}}  =O(1) ,\quad \abs{t}\geq 1.
\end{equation}
The equation (\ref{Zgtdef2}) and the bounds (\ref{HgtMuTbound3}) and (\ref{stirlinghalf}) together imply
\begin{equation}\label{ZgtuTbound2}
    \int_{\abs{t}\geq 1} \abs{Z_g(t)}u_t(t)\dt\ll_N T^{-N},
\end{equation}
and for $\abs{t}\leq 1$ we have
\begin{equation}\label{ZgtuTbound3}
    \int_{\abs{t}\leq 1} \abs{Z_g(t)}u_t(t)\dt=\int_{\abs{t}\leq 1} \abs{L_{p/q}\left(\frac{1}{2}+it\right)}u_T(t)\dt=O(1).
\end{equation}
The bound (\ref{ZgtuTbound2}) and (\ref{ZgtuTbound3}) together gives
\begin{equation}\label{ZgtuTbound4}
    \inR \abs{Z_g(t)}u_t(t)\dt=O(1)
\end{equation}
as $T\to\infty$. This contradicts the bound (\ref{ZgtuTbound1}).

\end{document}